\documentclass{ecos_2023}

\title{Towards CO2 valorization  in a multi remote renewable energy hub framework}

\author{V. Dachet\textsuperscript{a}, A. Benzerga\textsuperscript{b}, R. Fonteneau\textsuperscript{c}, D. Ernst\textsuperscript{d}\textsuperscript{e}}

\address{
\textsuperscript{a} University of Liège,
Liège,
Belgium,
victor.dachet@uliege.be
\textbf{CA}
\and
\textsuperscript{b} University of Liège,
Liège,
Belgium,
abenzerga@uliege.be
\and
\textsuperscript{c} University of Liège,
Liège,
Belgium,
raphael.fonteneau@uliege.be
\and
\textsuperscript{d} University of Liège,
Liège,
Belgium,
dernst@uliege.be
\and
\textsuperscript{e} Telecom Paris, Institut Polytechnique de Paris,
Paris,
France,
damien.ernst@telecom-paris.fr}


\usepackage[utf8]{inputenc}
\usepackage{amsmath}
\usepackage{amsfonts}
\usepackage{amssymb}
\usepackage{graphics}
\usepackage{graphicx}
\usepackage{url}
\usepackage{xcolor}

\usepackage{hyperref}

\usepackage[title]{appendix}

\usepackage[numbers]{natbib} 
\abstract{\normalsize
\begin{abstract}

In this paper, we propose a multi-RREH (Remote Renewable Energy Hub) based optimization framework. This framework allows a valorization of CO2 using carbon capture technologies. This valorization is grounded on the idea that CO2 gathered from the atmosphere or post combustion can be combined with hydrogen to produce synthetic methane. The hydrogen is obtained from water electrolysis using renewable energy (RE). Such renewable energy is generated in RREHs, which are locations where RE is cheap and abundant (e.g., solar PV in the Sahara Desert, or wind in Greenland). We instantiate our framework on a case study focusing on Belgium and 2 RREHs, and we conduct a techno-economic analysis. This analysis highlights, among others, the interest of capturing CO2 via Post Combustion Carbon Capture (PCCC) rather than only through Direct Air Capture (DAC) for methane synthesis in RREH. By doing so, a notable reduction of 9.2\% is observed in the total cost of the system under our reference scenario. In addition, we use our framework to derive a carbon price threshold above which carbon capture technologies may start playing a pivotal role in the decarbonation process of our industries. For example, this price threshold may give relevant information for calibrating the EU Emission Trading System so as to trigger the emergence of the multi-RREH.

\end{abstract}
}

\usepackage{siunitx}
\usepackage{eurosym}

\DeclareSIUnit{\EUR}{\text{\euro}}
\usepackage[normalem]{ulem} 

\def\equationautorefname#1#2\null{%
  #1(#2\null)
}

\begin{document}

\maketitle

\section{Introduction}
\label{sec:intro}

While the whole world is engaged in a process to decrease greenhouse gas emissions, capturing CO2 appears more and more as a crucial element to limit global warming. Once it is captured, CO2 may be either stored (CCS - Carbon Capture and Storage), either valorized (CCU - Carbon Capture and Utilisation), for instance through synthetic methane generation. In this article, we focus on CCU, where CO2 is seen as a required ingredient in the process of generating synthetic methane, together with \textit{green} hydrogen, i.e. hydrogen obtained from renewable energy-based electrolysis.

In this paper, we build on top of the Remote Renewable Energy Hub (RREH) approach  \cite{Berger2021} to propose a multi-hub, multi CO2 sources approach. CO2 is captured using both Post-Combustion Carbon Capture (PCCC) and Direct Air Capture (DAC) technologies. Hydrogen is produced from electrolysis using renewable energy in a RREH which is particularly well-suited for producing cheap and abundant renewable energy (e.g., solar energy in the Sahara desert, or wind energy in Greenland). The RREH concept also relies on the following idea: some locations show large amount of energy consumption while not having lots of renewable energy resources (e.g., Europe). On the opposite, some places have abundant renewable energy while having almost no energy demand. In its original formulation, the RREH concept suggests to use DAC technologies to feed the CO2 demand at the RREH. In this paper, we include PCCC technologies as an alternative to DAC technologies: in addition or replacement to being captured in the atmosphere, CO2 emitted in energy intensive locations may be transported to the RREHs to be combined with green hydrogen for producing neutral synthetic methane.

We propose a methodology for assessing the technico-economic feasibility of exporting CO2 into RREH where synthetic CO2-neutral methane would be generated using locally produced green H2. We formalise an optimisation problem where CO2 sources are in "competition" to provide CO2 to the methanation units in the RREHs. This methodology is based on a linear program modelling of Belgium energy system, including gas and electricity demand, and main CO2 emitters. We rely on previously published approaches to develop our approach Berger et al. \cite{Berger2021}, and, in particular, we use the GBOML language Miftari et al. \cite{Miftari2022} to model the energy system and to optimize it.

Our methodology is evaluated in the Belgian context: we consider Belgian CO2 emissions and Belgian gas and electricity demand. CO2 may be captured using Post Combustion Carbon Capture (PCCC) in Belgium or DAC in RREH locations. CO2 neutral synthetic methane will be produced in a remote energy hub from where it would be shipped back to serve the Belgian gas demand. We derive a CO2 emission cost in order to have a neutral emission system. We also determine a value of lost load (\textit{i.e.} a price associated with a lack of energy service) in order to serve the energy demand at all times. Several scenarios are studied with different prices of CO2 emissions, allocation or not of unserved energy and forcing of a given RREH.

\section{Related Work}
\label{sec:related_work}

This work is mainly related with the following topics that may play an important role in the deep decarbonation of our societies: (i) global grid  approaches, (ii) power-to-X technologies, multi-energy systems and and energy hub approaches, and (iii) CO2 quotas markets.

Global Grid (GG) approaches \cite{chatzivasileiadis2013global}, \cite{yu2019global}, sometimes referred to as Global Energy Interconnection approaches \cite{liu2015global},  are related with the idea of harvesting renewable energy from abundant and potentially remote renewable energy fields to feed the electricity demand in high demand centres. 
These approaches have mainly been oriented towards solutions using the electricity vector to repatriate energy from energy hubs, and have received a growing interest starting from the DESERTEC concept \cite{samus2013assessing} that focuses on Sahara solar energy resources from the Sahara desert to serve the European electricity demand. More recently, wind from Northern Europe and Greenland has also been identified as a promising resource to be valued within the GG context \cite{radu2022assessing}. Resource and demand configurations combining several types of resources as well as demand time zones show better results \cite{yu2019global}.

Multi-energy systems approaches \cite{munster2020sector,o2016energy} exploit the benefits of integrating energy demand and generation, as well as infrastructure. Power-to-X technologies, in particular power-to-CH4 technologies using hydrolysis and renewable energy for producing H2 \cite{GOTZ20161371}, offer a CO2 neutral solution to serve gas demand, but also a way to store vast quantities of energy issues from renewable sources \cite{BLANCO20181049}. Recently, Berger et al. have proposed a modeling framework \cite{Berger2021} for assessing the techno-economics viability of carbon-neutral synthetic fuel production from renewable electricity in remote areas where high-quality renewable resources are abundant. Let us mention that the idea of energy hubs was preexisting the work of Berger et al. \cite{geidl2006energy,mohammadi2017energy,sadeghi2019energy}, however the contribution of Berger et al. is the introduction of remote energy production, far from the demand. Our contribution is in line with the latter.


As this work aims to enhance the value of CO2, it is closely linked to the European Union Emissions Trading System (EU ETS). The EU ETS system, which is described on the European Commission's website \footnote{\url{https://climate.ec.europa.eu/eu-action/eu-emissions-trading-system-eu-ets_en}} and in \cite{carbonMarkets}, is a 'cap and trade' program. The system sets a cap on the total amount of certain greenhouse gases (GHG) that can be emitted by the facilities covered by the ETS. Within the cap, facilities are given emissions allowances, which can be traded with one another. The total number of allowances available is limited to ensure that they have value, and the cap is gradually reduced over time to lower total emissions. If a facility fails to cover its emissions fully, it faces substantial fines. Conversely, if a facility reduces its emissions, it can either retain the surplus allowance for future use or sell it to another facility that has not succeeded in covering its own emissions. This trading mechanism aims to reduce GHG emissions as soon as it becomes the most cost-effective solution and encourage investments in low GHG emissions solutions.

\section{CO2 Valorisation in a Multi-Remote Renewable Energy Hubs Approach}
\label{sec:Multi-RREH}

The Remote Renewable Energy Hub concept was first introduced in \cite{Berger2021} where the authors proposed a hub for synthesizing CH4 based on hydrogen and CO2 captured from the air thanks to a methanation unit. This concept has emerged within the context of global grid \cite{chatzivasileiadis2013global} and multi-energy systems approaches. These approaches aim at optimising the generation and utilisation of renewable energy (RE) by both (i) looking for abundant and cheap RE fields, (ii) taking advantage of daily/seasonal complementary of RE, as well as (iii) using power-to-gas technologies for better addressing RE generation fluctuations and meet e-fuels demand to act as a substitute for molecules derived nowadays from fossil fuels.


In the original article \cite{Berger2021}, the methanation unit was supplied with CO2 by a Direct Air Capture unit, and the energy demand was fulfilled by a single RREH located in Algeria. However, in this paper, we propose to investigate the feasibility of valorizing CO2 captured through Post Combustion Capture techniques at the energy demand center (EDC). Additionally, we deviate from the original paper by introducing a multi-RREH approach, wherein the EDC serves as a CO2 provider to a set of multiple RREHs, denoted as ${RREH_1, \ldots, RREH_h }$. Each hub $RREH_i ( 1 \leq i \leq h )$ has its unique characteristics, such as renewable energy type, potential, distance from the EDC, and means of CO2 transport from the EDC, which can affect its competitiveness.



\begin{figure}[p]
    \vspace{-3ex}
    \centering   \includegraphics[width=\textwidth, height=\textheight,keepaspectratio]{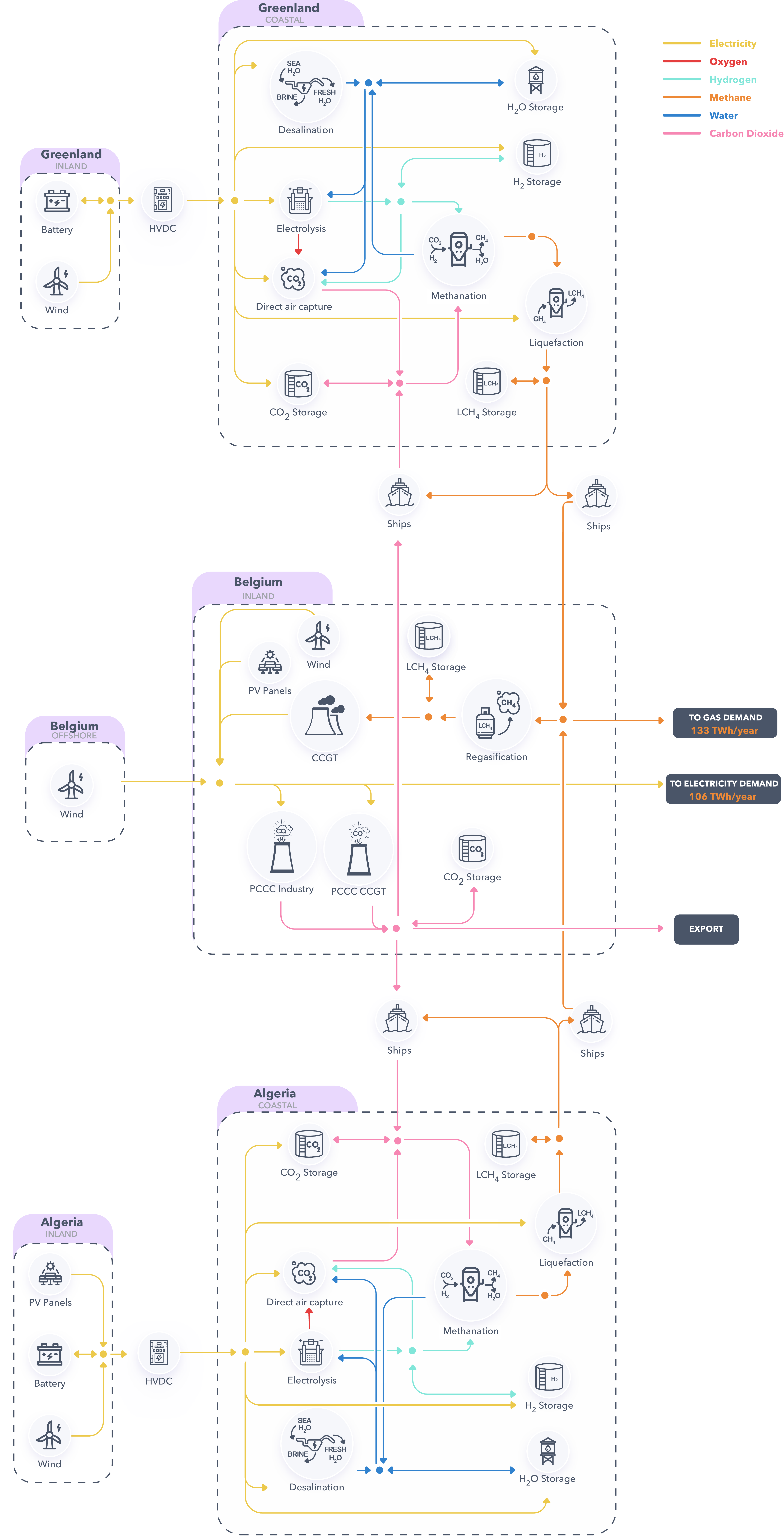}
    \caption{A schematic illustration of the remote energy hub. CO2 being captured, it may be used to synthesize fuel either locally either in a remote energy hub where renewable energy may be cheaper and more abundant.}
    \label{fig:RREH_model}
\end{figure}

In order to illustrate the concepts discussed above, we have developed a model for a multi-RREH system based on the following assumptions: (i) the EDC is Belgium, encompassing its gas and electricity demands as well as its CO2 emissions, (ii) there are two RREHs: one situated in the Sahara desert with access to solar and wind resources, and another in Greenland benefiting from the high-quality wind fields in the region. A detailed schematic of the resulting system is shown in \autoref{fig:RREH_model}. Similar to \cite{Berger2021}, we employed the GBOML language \cite{Miftari2022}, a recently developed language tailored for energy system optimization (refer to \autoref{sec:modelling} for more information), to model the system.

We note that the GBOML model code with two RREHs and one EDC system is available online\footnote{\url{https://gitlab.uliege.be/smart_grids/public/gboml/-/tree/master/examples}} and can  be easily extended to add additional RREHs and EDCs.

\section{Modelling}\label{sec:modelling}

This section provides insight into the optimization framework that underlies the multi-energy system model proposed in this work. The GBOML language introduced in \cite{Miftari2022}, a recently developed language dedicated to modeling graph-based optimization of multi-energy systems, is utilized to build this model. The optimization problem can be viewed as optimization on graphs, where a multi-energy system is considered as a set of nodes $\mathcal{N}$ that contribute to the (linear) objective and local constraints, and hyperedges $\mathcal{E}$ are used to model the constraints between nodes, such as those between RREHs and the EDC in our context.

The formalism employed in this work follows that introduced in \cite{Berger2021}. The entire system is defined by sets of nodes $\mathcal{N}$ and hyperedges $\mathcal{E}$. The optimization horizon is denoted by $T$, with time-steps indexed by $t \in \mathcal{T}$, where $\mathcal{T} = \{1, \ldots, T\}$.

A node $n \in \mathcal{N}$ is defined by internal $X^{n}$ and external $Z^{n}$ variables, where internal variables describe the specific characteristics of the unit, such as the nominal power capacity installed in the asset. Equality constraints $h_i(X^{n}, Z^{n}, t)=0$ with $i \in \mathcal{I}$ and inequality constraints $g_j(X^{n}, Z^{n}, t) \le 0$ with $j \in \mathcal{J}$, are employed for each $t \in \mathcal{T}$ to model operational constraints.

Each node $n$ has an associated cost function $F^{n}(X^{n}, Z^{n}) = \sum_{t=1}^{T} f^{n}(X^{n},Z^{n},t)$ that typically represents the capital expenditure and operational expenditure, i.e., CAPEX and OPEX, respectively.

Finally, equality and inequality constraints on hyperedges can be defined as $H^{e}(Z^{e}) = 0$ and $G^{e}(Z^{e}) \le 0$ with $e \in \mathcal{E}$ to model the laws of conservation and caps on given commodities.

    One can read this type of problem as:
     \begin{equation}\label{eq:prob_statement}
    \begin{aligned}
    \min \quad &  \sum_{n=1}^{N} F^{n}(X^{n}, Z^{n})\\
    \textrm{s.t.} \quad & h_i(X^{n}, Z^{n},t) = 0, \forall n \in \mathcal{N}, \forall t \in \mathcal{T}, \forall i \in \mathcal{I}\\
      \quad & g_j(X^{n},Z^{n},t) \le 0, \forall n \in \mathcal{N}, \forall t \in \mathcal{T}, \forall j \in \mathcal{J}\\
      \quad & H^{e}(Z^{e}) = 0, \forall e \in \mathcal{E}\\
      \quad & G^{e}(Z^{e}) \le 0, \forall e \in \mathcal{E}. \\
    \end{aligned}
    \end{equation}


The main assumptions underlying our model are the following:
\begin{itemize}
    \item Centralised planning and operation: In this framework, a single entity is responsible for making all investment and operation decisions.
    \item Perfect forecast and knowledge: It is assumed that the demand curves, as well as weather time series, are available and known \textit{in advance} for the entire optimisation horizon, i.e., $\forall t \in \{ 1, \ldots , T \}$.
    \item Permanence of investment decisions: Investment decisions result in the sizing of installation capacities at the beginning of the time horizon. Capacities remain fixed throughout the entire optimisation period, i.e., $\forall t \in \{ 1, \ldots , T \}$.
    \item Linear modelling of technologies: All technologies and their interactions are modelled using linear equations within this framework.
    \item Spatial aggregation: The energy demands and generation at each node are represented by single points. The topology of the embedded network required to serve this demand locally is not modelled in this approach. This can be viewed as an extension of the copper plate modelling approach used in electrical power systems.
\end{itemize}

     In our problem, all cost functions and constraints are affine transformation of the inputs. More details on the constraints of each technology can be found in \cite{BERGER_power_to_gas_2020106039}, \cite{Berger2021}.
    Additionnaly, the local objective function corresponding to the CAPEX is modelled with a uniform weighted average cost of capital (WACC) of $7\%$ for each technology. Thus, the CAPEX is computed using the following formula:
    \begin{equation}
        \zeta^n = \mbox{CAPEX}_n \times \frac{\mbox{w}}{(1 - (1 + \mbox{w})^{-\mbox{L}_n})}
        \label{waccformula}
    \end{equation}
    with $L_n$ the lifetime of technology n and w the WACC. Hence, $\zeta^n$ represents the annualised cost of investing in technology n.

    Moreover, a cap on the net CO2 emissions (\textit{i.e.} release in minus captured from the atmosphere) is added to the model. This latter is defined as 
    \begin{equation}\label{eq:cap_co2}
        \sum_{t \in \mathcal{T}} ( \sum_{a \in \mathcal{A}} q_{co2, t}^{a} - \sum_{c \in \mathcal{C}} q_{co2, t}^{c} ) \le \kappa_{co2} \nu
    \end{equation}
    with $\mathcal{A}$ and $\mathcal{C}$ representing the sets of technologies that release CO2 into the atmosphere and those that capture CO2 directly from the atmosphere, respectively, $\kappa_{co2}$ represents the CO2 cap in kilotons per year, and $\nu$ represents the number of years covered by the optimization horizon. The shadow price, or marginal cost, which is the dual variable associated with \autoref{eq:cap_co2} allows for the derivation of a CO2 cost in €/t. A detailed explanation of dual variables as marginal costs in linear programming can be found in \cite[Chapter 4]{linearOptim}.

\section{Case Study: Belgium}

This case study is focused on Belgium with two remote renewable energy hubs: one located in Algeria and another one located in Greenland. We will analyse the techno-economic feasibilty of the system while responding to an energy demand composed only of electricity and gas in Belgium.

\vspace{1ex}
\subsection{Data}
The data cover 2 years: 2015 and 2016. It is used to characterize energy demand as well as load factors for renewable energy sources. 

\textbf{Renewable generation profiles}

In order to determine the generation profiles of variable energy sources in Belgium we use the data from the transmission system operator (TSO) of Belgium \cite{Elia_power}. The profiles for the RREH located in Algeria are extracted with the same methodology as in \cite{Berger2021}. For the RREH situated in Greenland, the profiles of renewable energy are extracted thanks to the MAR model \cite{MAR} and given a power curve for an offsore wind turbine MHI Vestas Offshore V164-9.5MW.

\textbf{Energy consumption}

\begin{figure}[h]
    \centering
    \includegraphics[width=0.95\textwidth]{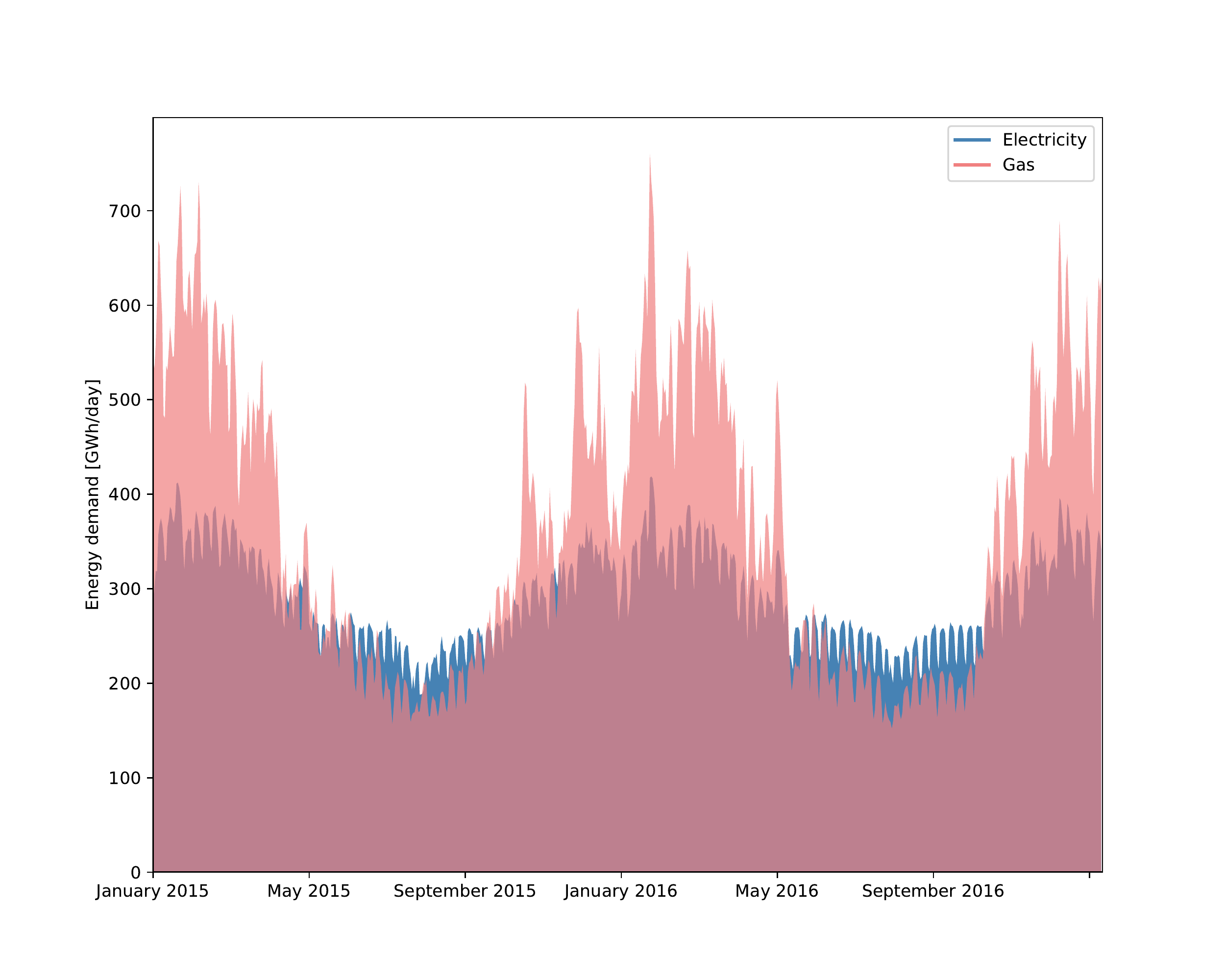}
    \caption{Daily aggregated profiles of electricity and natural gas demand covering the years 2015 and 2016 spanned by the optimisation.}
    \label{fig:energy_demand}
\end{figure}

The energy consumption data is collected for two energy vectors: gas (\cite{Fluxys}) and electricity (\cite{Elia_load}) with the same methodology as in \cite{BERGER_power_to_gas_2020106039}. In \autoref{fig:energy_demand}, the data corresponding to the two years is represented, where the signal is daily aggregated. In some cases, gas usage is shifted towards electricity needs, as described in \cite[section 4.2.2]{BERGER_power_to_gas_2020106039}. This shift is due to the use of heat pumps, which can help decarbonize heating in Europe. For both energy vectors, industrial and heating demands are taken into account. 

The peak power demand is equal to 60.13 GWh/h for both gas and electricity. The energy demand for electricity ranges from 6.42 to 20.29 GWh/h, while that for gas ranges from 5.51 to 39.84 GWh/h. The total energy demand is on average 106.45 TWh/year and 132.65 TWh/year for electricity and gas, respectively.

\vspace{1ex}
\subsection{Model Configuration}\label{subsec:model_config}

Our model consists of three main components (see \autoref{fig:RREH_model}): the energy demand center located in Belgium and two Remote Renewable Energy Hubs (RREHs) situated in Algeria and Greenland. The RREH in Algeria is modeled as described in \cite{Berger2021} with the same techno-economic parameters. The distinction is made with the inclusion of the CO2 connection between Belgium and Algeria. The RREH in Greenland is similarly modeled, with the exception of the removal of the photovoltaic potential and the modification of the high-voltage direct current (HVDC) line to a length of 100 km rather than 1000 km. 

The transportation of CO2 is achieved through the use of boats, which have a CAPEX of 5M€/kt, a lifespan of 40 years, and an average daily energy consumption of 0.0150 GWh/day. CO2 transport data was obtained from \cite{DanishEnergyAgency}. The loading and traveling time for these boats are assumed identical to those for liquefied methane carriers \cite{Berger2021}, \textit{i.e.} 24 and 116 hours, respectively. In order to fill the tank of CO2 carriers with fuel (liquefied methane), these tanks are loaded when unloading the CO2 at the RREH. Indeed, at the RREH, synthetic CH4 is available without having undergone any additional transport-related losses. Except for the storage facilities, liquefaction of CO2 has been excluded from the model. Sideways analyses have confirmed that this assumption has a negligible impact on the optimal objective.

Belgium is modeled with an electricity and gas demand as depicted in \autoref{fig:energy_demand}, with various means of production, including wind power, solar power, and a combined cycle gas turbine. The solar potential is limited to 40GW. The wind potential is equal to 8.4 GW and 8 GW for onshore and offshore capacities, respectively. The techno-economic parameters of each technology deployed in Belgium follow those in \cite{BERGER_power_to_gas_2020106039}.

We have also added a CO2 source that is equivalent to 40Mt CO2/year, which corresponds to the energy sectors and industrial processes greenhouse gases in Belgium in 2019 \cite[Table 4.1.1 (pp. 165- 166)]{EU-commission}. We assume that we can install post-carbon capture technologies (PCCC) in these sectors.

In terms of carbon capture technologies, the model has access to direct air capture installed at the RREHs, as well as a PCCC in Belgium on the 40Mt of CO2 per year and a PCCC installation on the CCGT.

As stated in \cite{BERGER_power_to_gas_2020106039}, the cost of PCCC is 3150M€/kt/h of CAPEX. The variable operating and maintenance costs (VOM and FOM) have been neglected in this analysis. However, a demand of $0.4125 GWh_{el}/kt_{CO2}$ of electricity is required. The expected lifetime is assumed to be 20 years.

Similarly, according to \cite{Berger2021}, the cost of DAC is equal to 4801.4 M€/kt/h of CAPEX. Similar to PCCC, VOM and FOM are ignored. The operational requirements for DAC are $0.1091 GWh_{el}/kt_{CO2}$ of electricity, $0.0438 kt_{H2}/kt_{CO2}$ of di-hydrogen, and $5.0 kt_{H20}/kt_{CO2}$ of water. The expected lifetime is assumed to be 30 years.

\vspace{1ex}
\subsection{Results}

In this section, we explore several scenari. We describe the variables that are used to differentiate the scenari
\begin{enumerate}
    \item Cost or Cap on CO2: either a cap is set of 0 t/year or a price at 80€/t or 0€/t
    \item Cost of energy not served (ENS): either energy not served is not allowed or a penalty of 3000€/MWh is imposed for each unit of unproduced energy.
    \item Forcing or not the use of a given RREH.
\end{enumerate}

\newpage
The results are generated with 5 scenari:

\textbf{Scenario 1}: This scenario seeks to avoid energy scarcity, whatever the cost. Therefore, no ENS is allowed. In addition, a hard constraint is set on CO2 emissions: a cap on CO2 is set.

\textbf{Scenario 2}:  This scenario follows the same assumptions as scenario 1 except that it leverages the constraint on energy not served. The cost associated to electricity not served is equal to 3000€/MWh, which is a standard value in the electricity context \cite{Voll}.

\textbf{Scenario 3}: This scenario leverages the constraint on CO2 emissions, and does not force the avoidance of energy not served but is penalized by 3000€/MWh not served. A penalty is associated with any CO2 emission in the atmosphere in the form of a fee equal to 80€/t - a value that reflects the current price of CO2 in the EU-ETS trading system \cite{CO2_price}.

\textbf{Scenario 4}: This scenario follows the same assumptions as scenario 3, with the difference that the cost of CO2 is equal to 0€/MWh. The aim is to showcase the system's configuration in the absence of any considerations for CO2 emissions.

\textbf{Scenario 5}: This scenario follows the same assumptions as scenario 1, with the difference that the only available RREH is in Greenland. 

\begin{table}[t] 
    \centering
    \begin{tabular}{|c|c|c|c|c|c|}
            \hline
            Scenario & Cap on CO2  & Cost of CO2  & ENS & Cost ENS  & Objective \\ 
              & (kt) & (M€/kt) &  & (k€/MWh) & (M€) \\
            \hline
            1 & 0.0 & 0.0 & No & - & 83742.61 \\
            2 & 0.0 & 0.0 & Yes & 3.0 & 80778.02 \\
            3 & No & 0.08 & Yes & 3.0 & 78872.94 \\
            4 & No & 0.0 & Yes & 3.0 & 76323.94 \\
            5 & 0.0 & 0.0 & No & - & 111209.95 \\
            \hline
            
    \end{tabular}
    \caption{Scenari parameters.}
    \label{tab:scenario_parameters}
\end{table}

These scenari summarized in \autoref{tab:scenario_parameters} vary in their degree of constraint. Scenario 1 is the most restrictive, with a cap on CO2 emissions and no allowance for energy not served. Scenario 2 allows for energy not served, while scenarios 3 and 4 remove the cap and replace it with CO2 prices of 80€ and 0€ per ton, respectively. Finally, scenario 5 requires the use of the RREH in Greenland, with parameters identical to those of scenario 1.

\vspace{1cm}
\begin{figure}[h]
    \centering
    \includegraphics[width=\textwidth]{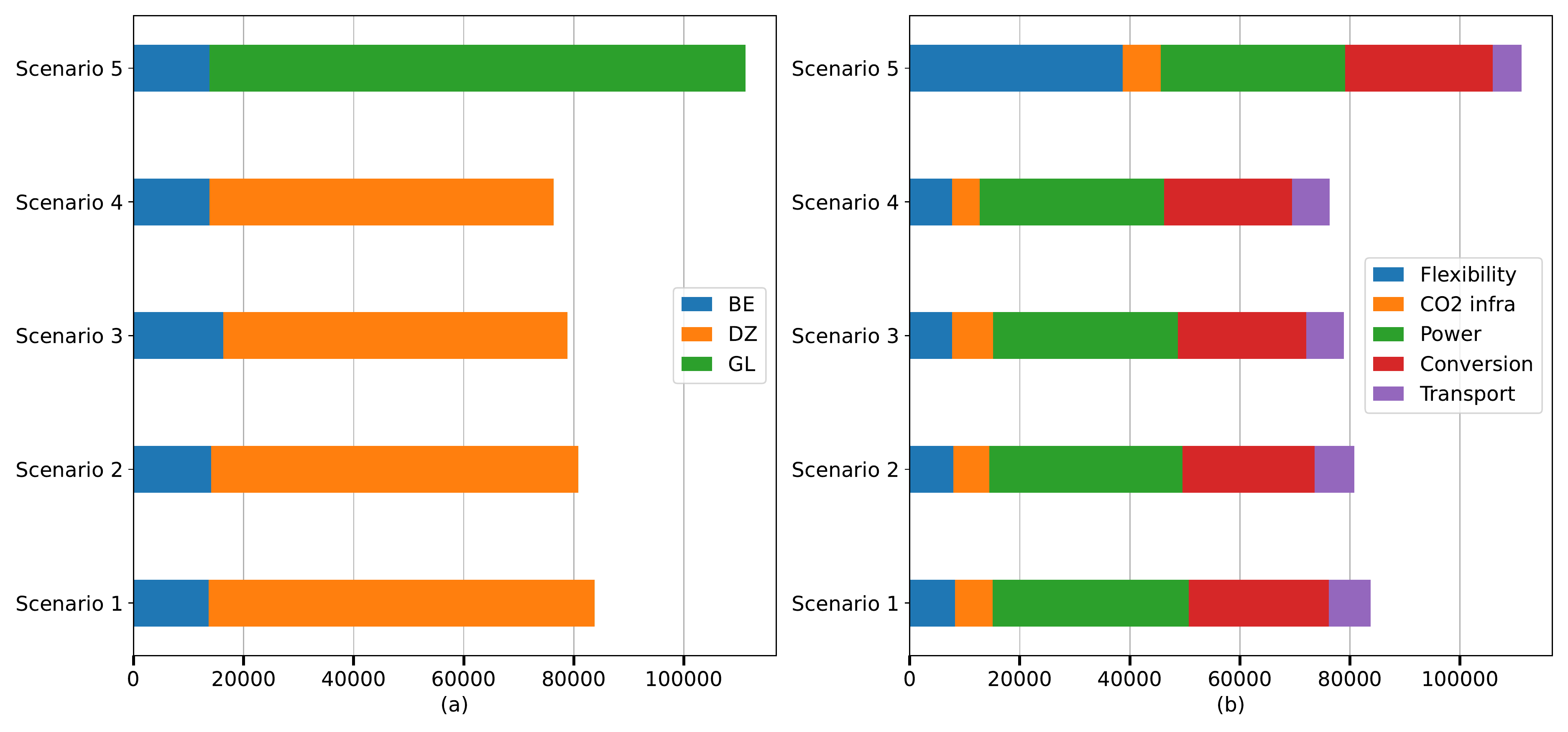}
    \caption{(a): Breakdown of costs per scenario and per cluster (Belgium (BE), Algeria (DZ), and Greenland (GL)). (b): Breakdown of costs per scenario per asset function. Flexibility covers storage capacities, CO2 Infra covers CO2 capture, storage, and transport, Power covers means of electricity production, Conversion covers all assets that convert one commodity into another and Transport HVDC lines and CH4 carriers.}
    \label{fig:cost_per_scenario_cluster}
\end{figure}

\vspace{1ex}
\subsection{Analyses and Discussion}
In this section, we introduce and discuss the results in detail. We choose to present a cross-scenario analysis in the light of key indicators and statistics extracted from the model.

\textbf{Total cost.}

The results indicate that the costs associated with enabling the hub in Algeria are substantially lower than those in Greenland, as depicted in  \autoref{fig:cost_per_scenario_cluster} (a) where nothing is built in the Greenland hub from scenarios 1 to 4, despite it being available for use. This disparity in costs can be attributed to the over-dimensioning of flexibility assets, particularly the storage capacities, as illustrated in \autoref{fig:cost_per_scenario_cluster} (b). This is mainly explained to electricity generated solely through wind available in Greenland, whereas both solar and wind electricity are obtainable in Algeria. This implies that the flexibility assets have to take the lead in maintaining a minimum of electricity delivery required in the electrolysis power plant. 

Furthermore, a reduction in total costs is observed in the first four scenarios with respect to the objective. This is explained with the order on the scenari based on their degree of constraint with scenario 1 being the most constrained and scenario 4 being the least.

\textbf{Power installation capacities.}

All power capacities installations are displayed in \autoref{tab:power}.

The potential in Belgium of solar energy is never reached while for both wind offshore and onshore the potential is reached in all scenari. 

From scenario 1 to scenario 2, the only difference being the allowance of ENS, there is an increase in the installation of controllable energy production assets. Indeed, there is a shift in capacity from CCGT to solar energy in Belgium between the first scenario and the second.

Comparing scenario 1 and 5, solar energy in Belgium is more expensive than importing CH4 from the RREH in Algeria. Importing from Greenland is more expensive and leads to an increase in power capacity installation in Belgium for solar, but it does not reach the maximum potential.

Another interesting comparison can be made with the work of \cite{Berger2021}, where the capacity installation in the hub for the reference scenario is 4.3GW of solar and 4.4GW of wind. In our case, the reference scenario 1 displays 100.51GW and 103.62GW, respectively. The power installation capacity is multiplied by approximately 23 while providing, on average, 282TWh/year of gas (HHV) to serve the gas demand and part of the electricity demand in Belgium, which is 28.2 times the gas production in the original paper.

\begin{table}[h]
    \centering
    \begin{tabular}{|c|c|c|c|c|c|c|c|}
         \hline
            Scenario & Wind onshore & Wind offshore & Solar & CCGT & Wind & Wind & Solar \\
             & BE & BE & BE & BE & GL & DZ & DZ \\\hline
            1 & 8.40 & 8.00 & 10.56 & 22.69 & 0.00 & 103.62 & 100.51 \\ 
            2 & 8.40 & 8.00 & 15.35 & 17.95 & 0.00 & 98.43 & 95.47 \\ 
            3 & 8.40 & 8.00 & 14.95 & 17.83 & 0.00 & 93.32 & 90.32 \\ 
            4 & 8.40 & 8.00 & 14.72 & 17.82 & 0.00 & 93.28 & 90.28 \\ 
            5 & 8.40 & 8.00 & 17.48 & 19.58 & 129.43 & 0.00 & 0.00 \\ 
            \hline
            
    \end{tabular}
    \caption{Total Power installation in GW per scenario.}
    \label{tab:power}
\end{table}

\textbf{CO2 installations (transport, capture).} 


In \autoref{tab:capture_co2}, the capacities of the CO2 capture units and the installations of transport capacity per scenario are displayed. Each time PCCC is activated, we recall that capturing CO2 is the only means to create gas in our system, and thus a minimum installation is required to support the demand. On the other hand, the DAC is only activated when a CO2 cap is set. PCCC has an efficiency of CO2 capture set to 90\%, which means that a direct air capture technology asset is necessary to recover the remaining 10\% of emissions in the atmosphere. This leads to a direct consequence, which is that when the DAC is available, the capacity of transport decreases because CO2 is locally available in the hub. However, the cost of CO2 capture by PCCC added to transport of CO2 is cheaper than the cost of DAC in the RREH. The only way to put PCCC out of business would be to have a distance between the hub and the energy demand center so long that the transport cost would increase too much.

Due to the higher concentration of CO2 in manufacturing smoke compared to the air, PCCC will likely always be cheaper than DAC, even with significant improvements in the DAC process. As a result, the operational costs associated with the energy required for PCCC will be lower than those of DAC.

\begin{table}[h]
    \centering
    \begin{tabular}{|c|c|c|c|c|c|c|c|}
            \hline
            Scenario & PCCC & PCCC CCGT & DAC DZ & DAC GL & Carrier DZ & Carrier GL \\ \hline
            1 & 4.11 & 2.62 & 1.30 & 0.00 & 8.030 & 0.000 \\
            2 & 4.11 & 2.07 & 1.47 & 0.00 & 7.142 & 0.000 \\
            3 & 4.11 & 1.80 & 0.00 & 0.00 & 9.694 & 0.000 \\
            4 & 3.76 & 2.06 & 0.00 & 0.00 & 9.701 & 0.000 \\
            5 & 4.11 & 2.40 & 0.00 & 1.35 & 0.000 & 7.564 \\
            \hline
            
    \end{tabular}
    \caption{Capacity, in kt/h, of CO2 capture technology and transport by hub and per scenario.}
    \label{tab:capture_co2}
\end{table}


\textbf{Cost of CO2 derived and Cap of CO2.}

From the first, second, and fifth scenarios, we are able to derive a shadow price thanks to the CO2 cap constraint. These correspond to approximately 162.77€/tCO2 for the first and second scenarios and 235.65€/tCO2 for the fifth scenario. This shows that given the system considered, i.e., Belgium and RREHs, putting a price of CO2 equal to 162.77€ would avoid these emissions in the atmosphere and activate the export of CO2 to Norway for storage purposes. In scenario 3, where a price of 80€/tCO2 is set, there is a net balance in the atmosphere of approximately 15Mt/year. In scenario 4, where no price is fixed, there is a net balance in the atmosphere which is equivalent to 16Mt/year. 

We would like to emphasize that the CO2 cap in our model only considers the emissions from the industrial and energy sectors, which are fully modeled. It does not account for a part of the emissions resulting from the gas demand served. Of this demand, 32\% is attributed to industrial needs, which are included in the statistics of the 40 Mt of CO2 emitted per year (see \autoref{subsec:model_config}), while the remaining 68\% is due to heating and is not covered by our cap. This heating gas demand translates to approximately 12.3 Mt of CO2 emitted per year.



\textbf{Cost of CH4 derived}

To estimate the cost of CH4 production, we first subtract from the optimal objective function the cost of the means of electricity production in Belgium (PV, on/off shore wind, CCGT), the cost of unserved energy (when applicable), and the cost related to export of CO2 for sequestration. All of these costs are substracted because they do not refer directly to the cost of producing synthetic methane. Then, we divide the obtained cost by the total energy content (HHV) in CH4 produced at the output of the regasification power plant in Belgium.

These methane costs, listed in \autoref{tab:ch4_price}, are compared to the price of 147.9€/MWh of methane (HHV) obtained by \cite{Berger2021}. Our scenarios achieve a lower cost for gas production (except for Greenland). This demonstrates that PCCC, which uses smoke with a high concentration of CO2 combined with transport, is more cost-effective than having only access to a DAC unit, as previously mentioned.

In our system, no fossil gas is available for import to Belgium; only synthetic gas produced from CO2 capture is used. If fossil gas were still available for import, our model would seek to minimize costs and import as much cheap gas as possible while staying within our carbon budget.

\begin{table}[h]
    \centering
    \begin{tabular}{|c|c|c|c|c|c|}
        \hline
        Scenario & 1 & 2 & 3 & 4 & 5 \\
        \hline
        [€/MWh] & 136.00 & 137.19 & 133.89 & 129.27 & 192.00 \\
        \hline
    \end{tabular}
    \caption{Estimation of methane price by retrieving the costs of power installations in Belgium, costs of unserved energy, and costs of exporting CO2 for storage purposes.}
    \label{tab:ch4_price}
\end{table}



\begin{figure}[h]
    \centering
    \includegraphics[width=0.47\textwidth]{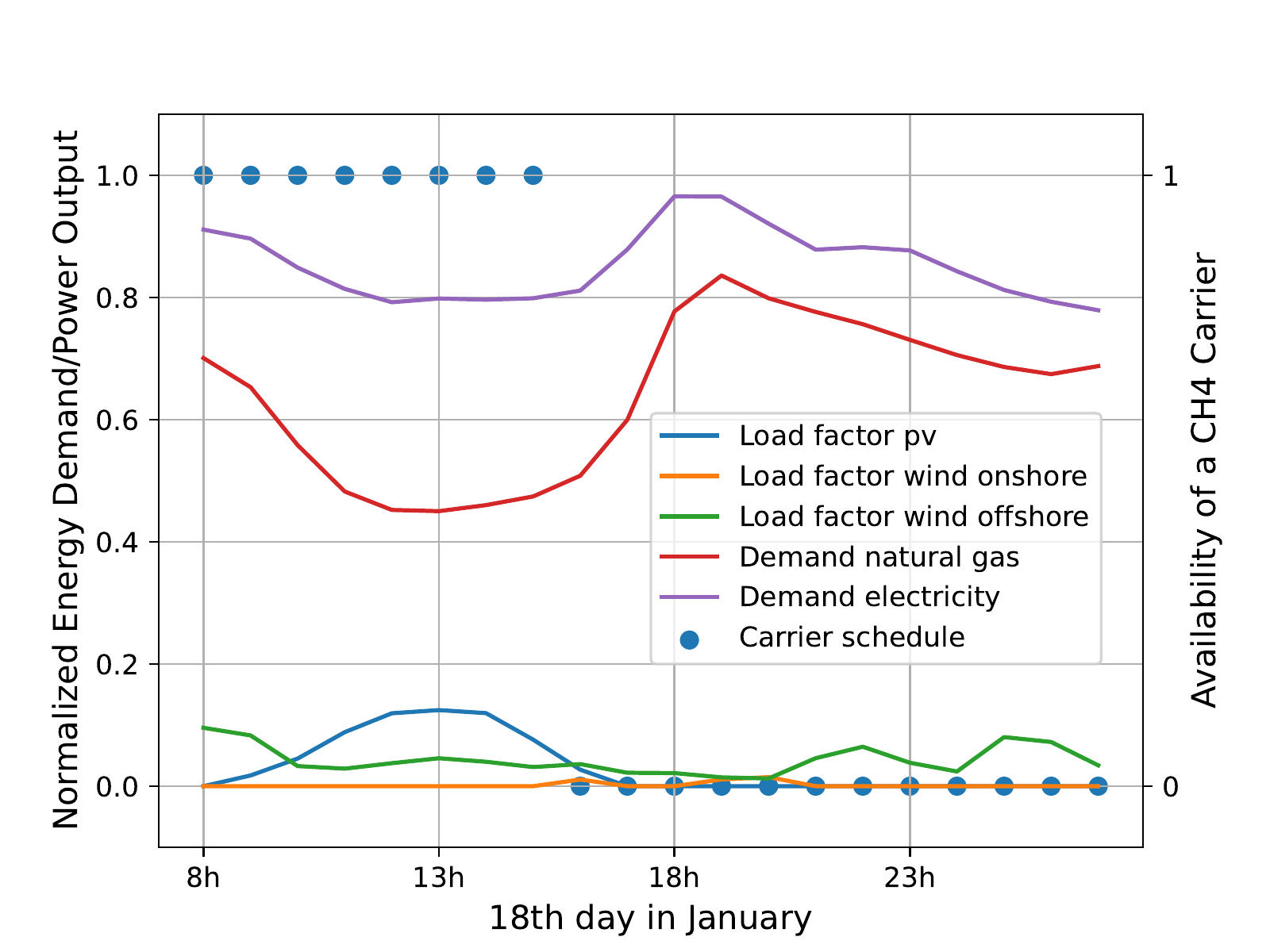}
    \caption{Evening of January 18th leading to the maximum shadow price associated with the hard constraint on energy not served in scenarios 1 and 5.}
    \label{fig:blackout}
\end{figure}

\textbf{ENS cost discussion}

The cost of unserved energy is a fixed parameter in scenarios 2, 3, and 4, but not in scenarios 1 and 5. Instead, a hard constraint is imposed to ensure that electricity demand is always met, resulting in a shadow price associated with the constraint. The maximum shadow price values for scenarios 1 and 5 are 913,640€/MWh and 1,075,913€/MWh, respectively. This is attributed to the peak in electricity and gas demand observed on January 18th at 18:00 (as shown in Figure \ref{fig:blackout}), where renewable energy load factors were low. Thus, all energy demand had to be supplied by the Combined Cycle Gas Turbine (CCGT) and gas resources.



\section{Conclusion}

In this work, we present our framework of multi remote energy hubs with capture of CO2 enabled in an energy demand center and its valorization by synthesizing methane in remote renewable energy hubs. We demonstrate the feasibility of serving the energy demand at the horizon 2050 of an entire country with only renewable energy and gas power plant fueled by synthetic methane while decarbonizing the energy and industry sectors on a case study implying Belgium as energy demand center and two RREHs: Greenland and Algeria. Our reference scenario exhibits a gas price of 136.0€/MWh instead of 149.7€/MWh in \cite{Berger2021} where only direct air capture was available in the RREH in order to feed CO2 into the methanation process. This shows the potential of Post Combustion Carbon Capture installations in the context of remote renewable energy hubs supply chains. We also derive a cost of CO2 of 163€ per ton in order to avoid any emission in the industrial and energy sector in Belgium. Finally, our model effectively captures the "competition" between different RREHs and is able to select exactly in which investments should be prioritized. In our simulations, the investments were made only for the RREH located in Algeria. In this respect, it would be interesting to study further how the different devices structuring the RREH in Greenland should be modified to become competitive with the RREH located in Algeria. This could be done for example by modifying the wind turbines selected for the Greenland hub so that they can operate with  higher nominal wind speeds and higher cut-off speeds in order to better exploit the strong winds in this area.

\section{Acknowledgements}
    The authors would like to thank Jocelyn Mbenoun for the templates and the useful conversations about energy, as well as Bardhyl Miftari and Guillaume Derval for their useful help with shadow pricing. The authors extend also their thanks to Julien Confetti for his precious help in the elaboration of script for generating the multi-hub picture. This research is supported by the public service of the Belgium federal government (SPF Économie, P.M.E., Classes moyennes et Energie) within the framework of the DRIVER project. Victor Dachet was supported by the Walloon region (Service Public de Wallonie Recherche, Belgium) under grant n°2010235 – ARIAC by
    \hyperlink{https://digitalwallonia4.ai/}{digitalwallonia4.ai}.

\vspace{2ex}

\begin{appendices}    
    \section{Glossary}
    \begin{tabbing}
        xxxxxxxxx\= xxxxxxxxxxxxxxxxxxxxxxxxxx \kill
        BE \> Belgium\\
        CAPEX \> Capital Expenditure\\
        CCGT \> Combined Cycle Gas Turbine\\
        DAC \> Direct Air Capture\\
        DZ \> Algeria\\
        EDC \> Energy Demand Center\\
        ENS \> Energy Not Served \\
        ETS \> Emission Trading System\\
        GBOML \> Graph Based Optimzation Modeling Language\\
        GL \> Greenland\\
        HHV \> Higher Heating Value\\
        OPEX \> Operational Expenditure\\
        PCCC \> Post Combustion Carbon Capture\\
        PV \> Photovoltaic \\
        RE \> Renewable Energy \\
        RREH \> Remote Renewable Energy Hub \\
        RES \> Renewable Energy Sources \\
    \end{tabbing}
\end{appendices}

\newpage
\section*{Nomenclature}

\textbf{Sets and indices}
\begin{description}[leftmargin=!,labelwidth=\widthof{maxlength}]
	\item[$\mathcal{E}, e$] set of hyperedges and hyperedge index
	\item[$\mathcal{G}$] hypergraph with node set $\mathcal{N}$ and hyperedge set $\mathcal{E}$
 
	\item[$\mathcal{I}^n, i$] set of external variables at node $n$, and variable index
 
    \item[$\mathcal{N}, n$] set of nodes and node index 
  
	\item[$\mathcal{T}, t$] set of time periods and time index

\end{description}

\textbf{Parameters}
\begin{description}[leftmargin=!,labelwidth=\widthof{maxlength}]
	\item[$\nu \in \mathbb{N}$] number of years spanned by optimisation horizon 
  
	\item[$\kappa_{i} \in \mathbb{R}_+$] maximum flow capacity of commodity $i$
  
    \item[$\zeta^{n} \in \mathbb{R}_+$] annualised CAPEX of node $n$ (flow component)
\end{description}

\textbf{Variables}
\begin{description}[leftmargin=!,labelwidth=\widthof{maxlength}]
    \item[$q_{it}^n \in \mathbb{R}_+$] flow variable $i$ of node $n$ at time $t$
\end{description}


  
  
  
  
  

  
  
        
        
   
  
  

\bibliography{bib.bib}

\end{document}